\documentclass{article}
\usepackage{iclr2025_conference_nopub,times}
\usepackage{amsmath,amssymb}
\usepackage{graphicx}
\usepackage{booktabs}
\usepackage{tabularx}
\usepackage{pifont}
\usepackage{hyperref}
\usepackage{url}
\usepackage{verbatim}

\usepackage{algpseudocode}

\usepackage{chngcntr}

\title{Continuous-Time Homeostatic Dynamics for Reentrant Inference Models}

\author{B. G. Chae \\
Electronics and Telecommunications Research Institute, Daejeon 34129, Republic of Korea \\
\texttt{bgchae@etri.re.kr}
}

\iclrfinalcopy
\begin{document}
\maketitle

\begin{abstract}
We formulate the Fast-Weights Homeostatic Reentry Network (FHRN) as a continuous-time neural-ODE system, revealing its role as a norm-regulated reentrant dynamical process. 
Starting from the discrete reentry rule $x_t = x_t^{(\mathrm{ex})} + \gamma\, W_r\, g(\|y_{t-1}\|)\, y_{t-1}$,
we derive the coupled system $\dot{y}=-y+f(W_ry;\,x,\,A)+g_{\mathrm{h}}(y)$
showing that the network couples fast associative memory with global radial homeostasis. 
The dynamics admit bounded attractors governed by an energy functional, yielding a ring-like manifold. 
A Jacobian spectral analysis identifies a \emph{reflective regime} in which reentry induces stable oscillatory trajectories rather than divergence or collapse. 
Unlike continuous-time recurrent neural networks or liquid neural networks, 
FHRN achieves stability through population-level gain modulation rather than fixed recurrence or neuron-local time adaptation. 
These results establish the reentry network as a distinct class of self-referential neural dynamics supporting recursive yet bounded computation.

\end{abstract}

\section{Introduction}
Biological neural systems operate as continuous dynamical processes shaped by recurrent excitation, homeostatic adaptation, 
and multi-timescale synaptic plasticity \citep{1,2,3,4}.
Cortical populations exhibit persistent attractor-like states, oscillatory flows, and normalization mechanisms that stabilize activity while preserving computational flexibility.
Such behavior has long been formalized within dynamical-systems frameworks, including Hopfield-type associative energy models \citep{5,6},
continuous recurrent neural ODEs \citep{7,8,9}, and population-level memory models \citep{10,11}, 
where feedback and homeostasis jointly prevent runaway amplification and enable stable working memory.

Beyond stability, cognitive neuroscience emphasizes the importance of \emph{reentrant signaling}:
reciprocal, self-referential exchange across distributed cortical pathways \citep{12,13}.
Under this view, recurrence is not merely a temporal buffer,
but a mechanism by which internal representations iteratively refine themselves—supporting reflection, imagination, and higher cognition.

In parallel, artificial recurrent frameworks explored related principles through fast-weight mechanisms. Early fast-weight models introduced rapidly updated
associative memory via outer-product plasticity rules \citep{14}, 
later extended to recurrent formulations \citep{15,16} and to fast-weight programmers enabled by linearized attention \citep{17}. 
More recent retention-based architectures—including Linear Transformers, RetNet, and self-referential weight matrices \citep{18,19,20}—provide efficient
temporal continuity but still lack explicit population-level reentrant feedback: 
internal activations do not re-enter the representational stream in a controlled, state-dependent manner, 
nor are they regulated by dynamical stability constraints.

A recent study introduced the Fast-Weights Homeostatic Reentry Network (FHRN) \citep{21}:
a discrete-time architecture combining low-rank fast memory, activity-dependent
normalization, and a learned reentry operator. 
Although that work empirically identified a stable ``reflective band'' under moderate feedback gain, its 
dynamical principles were characterized only through discrete metrics, 
leaving unanswered how stability, memory, and reentry interact as a continuous population flow.

In this work, we derive a continuous-time formulation of FHRN, 
showing that it realizes a two-timescale nonlinear neural dynamical system with explicit leakage,
reentrant excitation, and homeostatic radial stabilization. 
Through Jacobian spectral analysis and a Lyapunov energy argument, we prove bounded attractors and
establish input-to-state stability. These results situate the neural network within the broader landscape of neural dynamical systems—bridging biologically inspired
fast-weight mechanisms with continuous neural ODEs—and provide a physically interpretable foundation for stable self-referential computation.

\section{Theoretical Framework}

\subsection{Continuous-time formulation of reentry networks}
Figure 1 summarizes the discrete computational structure of the Fast-Weights Homeostatic Reentry Network. 
The network consists of three interacting components: (i) a fast associative memory trace $A_t$, (ii) a
state-dependent reentry operator $W_r g(\|y_t\|)$ acting on the current population
activity, and (iii) a nonlinear homeostatic term regulating radial activity norms.
Together these components define a recursive interpretation loop in which the latent
representation is repeatedly updated rather than merely propagated forward.

\begin{figure}[ht!]
\includegraphics[scale=0.55, trim= -3.0cm 15.5cm 0cm 0cm]{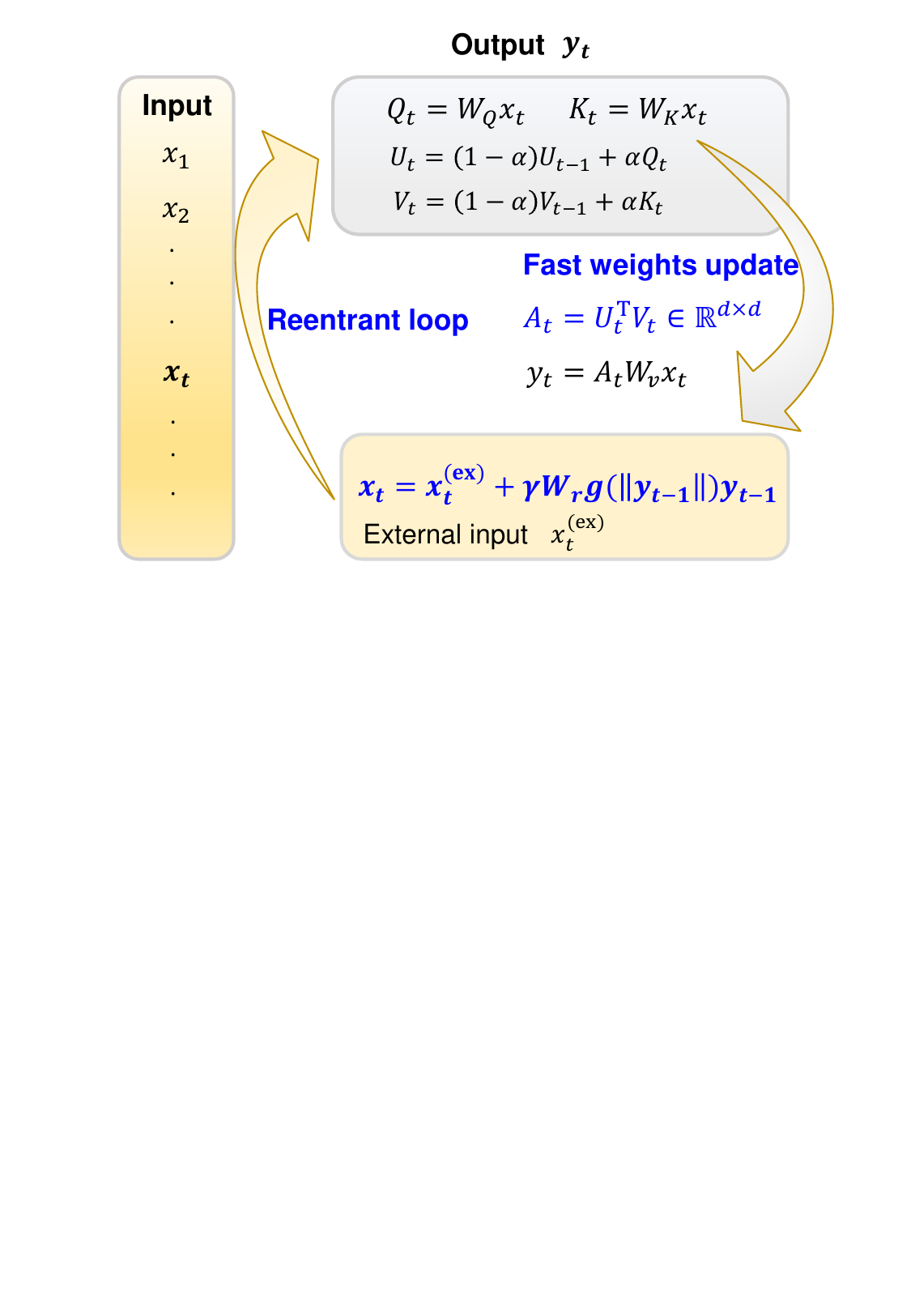}
\caption{
Computational architecture of the Fast-Weights Homeostatic Reentry Network.
At each timestep, the model receives an external token embedding $x_t^{(\mathrm{ex})}$ and combines it with a
reentrant signal $\gamma W_r g(\lVert y_{t-1}\rVert)\,y_{t-1}$ to form the effective block input $x_t$.
The queries and keys are derived via learned projections.
The fast associative memory $A_t = U_t^\top V_t$ stores low-rank temporal correlations
and modulates the internal representation $y_t$ on subsequent steps.
The reentrant loop provides recurrent reflective influence regulated by the homeostatic gain $g(\lVert y_t\rVert)$.
}
\end{figure}

\vspace{8pt}
The discrete-time FHRN update introduces reentrant feedback by injecting the
previous population activation back into the input representation through a
learned operator $W_r$:
\begin{equation}
x_t
= x_t^{(\mathrm{ex})}
+ \gamma\, W_r\, g(\|y_{t-1}\|)\, y_{t-1},
\label{eq:reentry_input}
\end{equation}
where $\gamma>0$ denotes the reentry gain. We define the population activity
radius as $r_t = \|y_t\|_2$,
and specify the homeostatic gain modulation as
\begin{equation}
g(r_t) = \frac{1}{1 + \beta (r_t^2 - 1)}, \qquad \beta > 0.
\label{eq:g_def}
\end{equation}
The term $g(r_t)$ decreases feedback strength when activity norms exceed the
target radius $r_t \approx 1$, providing a population-level normalization
mechanism that prevents runaway reentrant amplification.

When $\gamma = 0$, Eq.~\eqref{eq:reentry_input} reduces to a purely
feed-forward fast-weight transformer. For $\gamma>0$, the mapping becomes
self-referential: internal representations recursively re-enter the computation,
enabling iterative refinement rather than one-pass propagation.

Rather than implementing a classical hidden-state recurrence, this update defines a \emph{representational recurrence}: the external input $x_t^{(\mathrm{ex})}$ and the internally
generated feedback combine to form the instantaneous signal that is processed by the block mapping.
The output population activation is computed through a nonlinear transformation
\begin{equation}
    y_t \;=\; H(x_t),
\end{equation}
where $H(\cdot)$ includes fast-weight excitation, attention mechanisms, feedforward processing, normalization, and nonlinearity.  

\vspace{8pt}
During autoregressive inference, the external embedding evolves at a slow timescale (one token per step), whereas internal reentry operates on a faster implicit cycle.
Under this separation, $x_t^{(\mathrm{ex})}$ may be treated as locally constant.
Using a first-order multivariate Taylor expansion around 
$x_t^{(\mathrm{ex})}$ gives
\begin{equation}
    H(x_t)
    \;\approx\;
    H(x_t^{(\mathrm{ex})})
    \;+\;
    J_H(x_t^{(\mathrm{ex})})\,R(y),
    \label{eq:taylor_expansion}
\end{equation}
where $J_H(\cdot)$ denotes the Jacobian of the block transformation and $R(y)=\gamma\,W_r\,g(\|y\|)\,y$.
Defining the effective reentry operator
\begin{equation}
    W_r^\star := J_H(x_t^{(\mathrm{ex})})\,W_r,
\end{equation}
Eq.~\eqref{eq:taylor_expansion} reduces to the intrinsic reentrant form
\begin{equation}
    y \;\approx\; H(x_t^{(\mathrm{ex})}) \;+\; \gamma\,W_r^\star g(\|y\|)\,y,
\end{equation}
showing that the feedback loop acts directly in $y$-space rather than over raw inputs.

\vspace{8pt}
Writing the discrete update as a relaxation form, and taking the small-step limit,
yields the neural ODE
\begin{equation}
    \dot{y}
    =
    -y
    + \gamma\,W_r^\star y
    + g_{\mathrm{h}}(y)
    + c,
    \qquad
    c := H(x^{(\mathrm{ex})}),
    \label{eq:ode_with_drive}
\end{equation}
where the homeostatic correction acts radially and takes the form
\begin{equation}
g_{\mathrm{h}}(y)= -\kappa(\|y\|^2 - 1)\,y, \qquad \kappa>0.
\label{eq:g_homeo}
\end{equation}  
When analyzing local stability of the intrinsic reentrant flow,
the constant term $c$ representing the slow-timescale external drive may be omitted (\,$c=0$\,), yielding
$ \dot{y} \;=\; -y + \gamma\,W_r^\star y + g_{\mathrm{h}}(y)$,
which governs the autonomous internal dynamics of the neural network.

\vspace{4pt}
Collecting representational, fast-weight, and homeostatic contributions gives the final form
\begin{equation}
    \dot{y}(t)
    =
    -y(t)
    + f(Wy(t);\,x(t),A(t))
    + g_{\mathrm{h}}(y(t)),
    \label{eq:fhrl_final_ode}
\end{equation}
revealing that the network performs a continuous fixed-point search over its own reentrant
representation.  
The leak term ensures contraction, the reentrant fast-weight pathway provides structured excitation,
and the homeostatic field constrains trajectories to a bounded manifold—enabling stable and
self-referential dynamics.

\subsection{Jacobian spectral interpretation}
\vspace{4pt}
Consider the continuous dynamics in Eq. (9),
where $f(\cdot)$ encodes associative fast-weight excitation and 
$g_{\mathrm{h}}(y)$ provides activity-dependent radial damping.  
Linearizing around a stationary point $y_\star$ yields the Jacobian
\begin{equation}
    J_y(y) = -I 
    + F(Wy;x, A)\, W
    + D g_{\mathrm{h}}(y),
\end{equation}
with 
\[
F(Wy,x,A) 
= \left.\frac{\partial f}{\partial z}\right|_{z = Wy}
\]
denoting the gain of the $f$-path with respect to its preactivation.
The term $D g_{\mathrm{h}}(y)$ is the differential contribution of the homeostatic field,
\begin{equation}
    D g_{\mathrm{h}}(y) =-\kappa \left[(\|y\|^2-1) I + 2yy^\top \right],
\end{equation}
which supplies negative radial curvature in the supercritical region $\|y\|>1$.

Local stability requires all eigenvalues of $J_y(y_\star)$ to have negative
real parts.  Throughout, we use the notation
\[
\rho(J) \equiv \max_i \Re(\lambda_i(J)),
\]
i.e., the spectral abscissa of the Jacobian.  
Under this convention, the local stability condition becomes
\[
\rho(J_y(y_\star)) < 0.
\]

A sufficient condition can be written in small-gain form as
\begin{equation}
\|F(Wy_\star;x,A)\, W\|_2
<
1 - 
\lambda_{\min}(- D g_{\mathrm{h}}(y_\star)),
\end{equation}
meaning that the effective amplification through the $f$-path must be
smaller than the combined contraction produced by the intrinsic leak
$(-I)$ and the homeostatic dissipation.  
When this inequality holds, the reentrant flow is locally attracting and
trajectories converge toward a bounded invariant region.

If $f$ is Lipschitz in its preactivation, $\|\partial f/\partial z\|\le L_f$, then the same small-gain
condition ensures input-to-state stability:
\[
\|x(t)^{(\mathrm{ex})}\|_\infty < \infty 
\quad\Rightarrow\quad 
\sup_t \|y(t)\| < \infty.
\]
Thus bounded external inputs yield bounded neural trajectories.
Because the homeostatic field drives $\|y\|$ toward its preferred operating radius, even in the
presence of rotational components the flow remains trapped within a finite annulus of the state
space.  
These results show that the network admits well-defined stable reentrant dynamics under explicit,
interpretable conditions on the Jacobian spectrum.

\subsection{Lyapunov stability of reentrant dynamics}
To examine the stability of the continuous-time reentrant dynamics, we consider the candidate Lyapunov function
\begin{equation}
V(y) = \tfrac{1}{4}(\|y\|^2 - 1)^2,
\label{eq:lyapunov_def}
\end{equation}
which measures the deviation of the population activity from the homeostatic equilibrium radius.  
Differentiating along trajectories of continous-time ODE system yields
\begin{equation}
\dot{V} = (\|y\|^2 - 1)\,y^\top[-y + f(Wy;x,A) + g_{\mathrm{h}}(y)],
\label{eq:v_dot}
\end{equation}
where the direction of change remains aligned with the radial unit vector $\hat{y}=y/\|y\|$, since all corrective terms act radially.

The first (leak) term produces a strictly negative cubic dissipation term,
\begin{equation}
(\|y\|^2 - 1)\,y^\top(-y) = -(\|y\|^2 - 1)\|y\|^2,
\end{equation}
ensuring monotonic radial energy decay.
The nonlinear drive $f(Wy,x,A)$ captures bounded modulation from the fast-weight trace and external input; under layer normalization and limited feedback gain, it satisfies
\begin{equation}
\|f(Wy;x,A)\| \le c_A\|A\| + c_x\|x\|,
\label{eq:bounded_drive}
\end{equation}
for some constants $(c_A,c_x)$.  

The homeostatic term $g_{\mathrm{h}}(y)$ is precisely proportional to the gradient of the Lyapunov function.  
Near the equilibrium shell $\|y\|\simeq1$, it contributes an additional stabilizing term proportional to $(\|y\|^2-1)^2$.

Collecting these contributions, we obtain the differential inequality
\begin{equation}
\dot{V} \le -c(\|y\|^2 - 1)^2 + c_A\|A\|^2 + c_x\|x\|^2,
\label{eq:iss_bound}
\end{equation}
where $c>0$ encapsulates the joint effect of leak and homeostatic gains.  
This ensures that radial deviations from the equilibrium decay exponentially up to bounded perturbations due to $(A,x)$.  
In particular, for vanishing inputs $(A,x)\to0$, the state trajectories asymptotically approach the limit set $\|y\|=1$.  

Hence, the continuous reentrant dynamics are \textit{input-to-state stable} with respect to the external drive and associative memory perturbations.  
The homeostatic term enforces a Lyapunov-like global damping, while the reentry operator $W_r$ shapes the rotational component of the flow without affecting the boundedness of trajectories.

Physically, $V(y)$ acts as a smooth quartic energy basin centered on the homeostatic manifold $\|y\|=1$.  
The leak term $(-y)$ enforces dissipation, while $W_r$ induces structured rotational motion on this manifold.  
The resulting dynamics resemble a damped oscillator constrained to a curved activation shell, providing a continuous-time analog of stable reentrant inference.

\begin{figure}[ht!]
\includegraphics[scale=0.62, trim= -1.6cm 11.8cm 0cm 0cm]{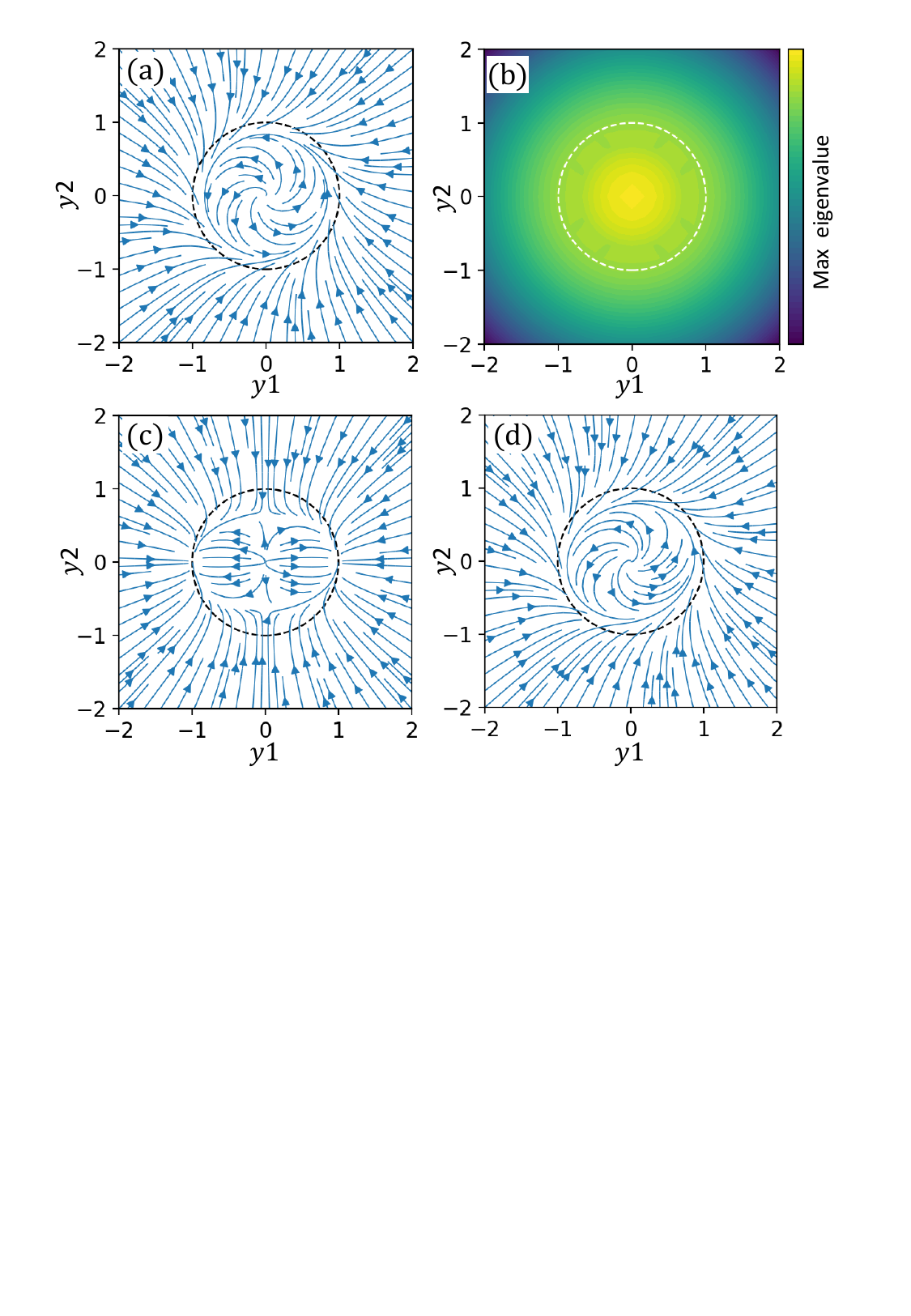}
\caption{
Phase–space behavior of the continuous FHRN dynamics under different forms of the reentrant operator.
(a) A purely antisymmetric operator yields rotational flow and a stable ring attractor.
(b) Local Jacobian spectral map showing the zero-real-eigenvalue contour aligning with the attractor radius.
(c) A symmetric operator produces saddle-like expansion and contraction without rotational structure.
(d) A mixed operator generates hybrid spiral dynamics, combining rotation and anisotropic gain.
}
\end{figure}

\section{Numerical analysis}

To analyze how the reentrant matrix shapes the flow, 
$W_r^\star$ is decomposed into its antisymmetric and symmetric parts,
\[
W_r^\star = W_{\!a} + W_{\!s}, \qquad 
W_{\!a} = \tfrac{1}{2}(W_r^\star - W_r^{\star\!\top}),\;
W_{\!s} = \tfrac{1}{2}(W_r^\star + W_r^{\star\!\top}).
\]
The antisymmetric component $W_a$ produces rotational flow conserving $\|y\|$, 
while the symmetric component $W_s$ generates stretching or compression along the principal axes.
Their superposition yields a hybrid spiral combining rotation and contraction.

Figure 2 illustrates how the reentrant operator $W_r^\star$ shapes the continuous-time FHRN dynamics by comparing phase flows generated 
by its antisymmetric, symmetric, and mixed components. 
In Fig. 2(a), a purely antisymmetric $W_r^\star$ induces a rotational flow characteristic of conservative dynamics. 
Blue streamlines show the temporal evolution of the population state in the two-dimensional activation space $y = (y_1, y_2)$, 
and the dashed circle denotes the homeostatic manifold $\|y\| = 1$.
Leakage $(-y)$ generates radial contraction, the reentrant term $W_r^\star y$ contributes angular momentum, 
and the nonlinear homeostatic field $g_{\mathrm{h}}(y)$ stabilizes the amplitude.
Together, these forces produce a contracting spiral that converges onto a stable ring attractor near $\|y\|\approx 1$, enabling persistent yet bounded rotational reentry. 
This regime mirrors cortical recurrent circuitry that maintains activity without runaway excitation.

Figure 2(b) shows the corresponding local Jacobian stability map, where color
encodes the maximum real eigenvalue of $J_y(y)$. The homeostatic shell aligns
with the zero-eigenvalue contour, confirming that the ring manifold forms a
structurally stable invariant set rather than a numerical artifact.

A symmetric reentrant operator, shown in Fig. 1(c), produces saddle-like
dynamics with expansion and compression along principal axes but with no
rotational component. In this regime, reentry becomes geometrically unstable
and cannot support sustained internal dynamics.
Finally, the mixed case in Fig. 1(d) exhibits a hybrid spiral: rotational motion
from the antisymmetric component combined with anisotropic expansion generated
by the symmetric one. Despite these competing effects, homeostasis restores
boundedness, steering trajectories toward the equilibrium ring.

Figure 3 illustrates the Lyapunov geometry and radial contraction properties 
that govern the stability of the FHRN continuous-time dynamics.
The Lyapunov potential $V(y)=\tfrac14(\|y\|^2-1)^2$ defines a Mexican-hat energy landscape with a minimum along the homeostatic ring $\|y\| = 1$, as illustated in Fig. 3(a).
The color surface encodes the energy magnitude, and the inserted box illustrates a representative trajectory descending the potential toward the stable orbit.
Unlike point-stable attractors, the minimum forms a one-dimensional invariant ring manifold, indicating that the dynamics converge radially while preserving tangential (rotational) degrees of freedom.

Figure 3(b) shows the temporal evolution of the activity norm $\|y(t)\|$
for different initial magnitudes.
States starting inside the shell ($\|y_0\|<1$) accelerate outward, while suprathreshold states contract inward, reflecting the radius-dependent curvature of the Lyapunov well.
Convergence is not purely exponential: trajectories exhibit fast decay near large radius and slower approach as they enter the shallow minimum near 
$\|y\|\approx 1$.
Together, the two plots confirm that the homeostatic field imposes strict radial stability while leaving angular dynamics—and hence reentrant structure—intact.

\begin{figure}[ht!]
\includegraphics[scale=0.78, trim= 0.5cm 19cm 0cm 0.5cm]{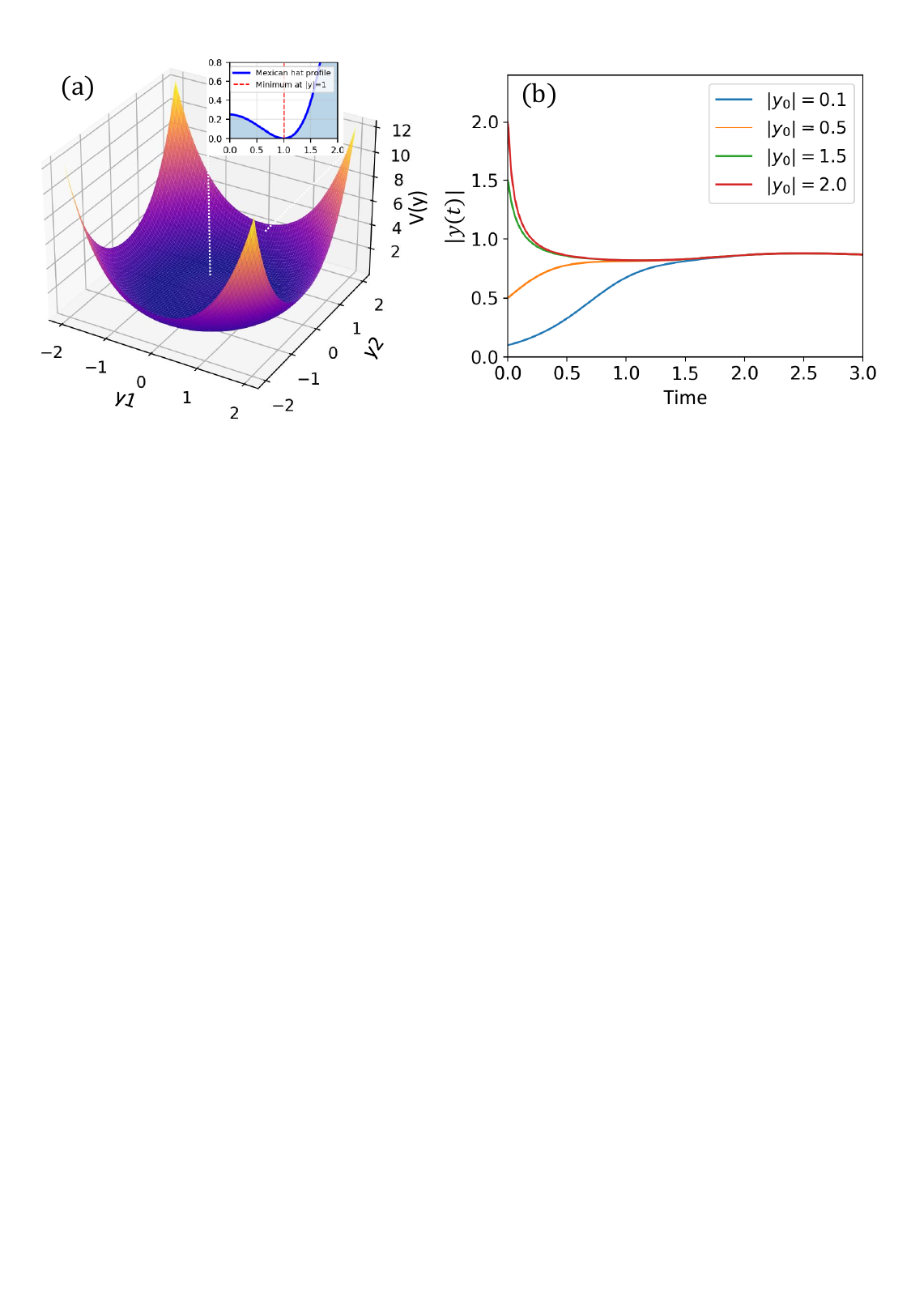}
\caption{
Stability of the continuous-time FHRN dynamics under the homeostatic field.
(a) Lyapunov energy landscape of the dynamics defined by \( V(y) = \tfrac{1}{4}(\|y\|^2-1)^2 \).
The system monotonically descends this energy surface toward the minimum manifold.
(b) Radial evolution $\lVert y(t)\rVert$ for different initial conditions, demonstrating monotonic regulation toward the homeostatic manifold. 
}

\end{figure}

Figure 4 illustrates the stability structure of the FHRN dynamics across the
parameter space defined by the homeostatic gain $\beta$, the reentry gain
$\gamma$, and the spectral magnitude of the reentrant operator $\|W_r\|$.
The 3D stability surface shows the critical boundary
$\|W_r\|_{\mathrm{crit}} = \frac{1}{\gamma\,g(\|y\|)}$ separating stable and unstable
dynamical regimes, as depicted in Fig. 4(a). 
Increasing either homeostasis ($\beta$) or the reentry gain
($\gamma$) lowers the allowable reentrant strength before instability emerges,
revealing an inverse trade-off between normalization and feedback amplification.

The corresponding 2D stability map in Fig.~4(b) provides a slice of this
surface at fixed radius ($\|y\| = 1.1$).  
To quantify the strength of the reentrant feedback, we plot the
\emph{effective reentry gain}
\[
r_e = \gamma\, g_{\mathrm{h}}(\|y\|;\beta)\, \|W_r\|_2,
\]
which captures when the reentrant drive becomes sufficiently strong to
counteract the intrinsic leak in the Jacobian 
$J = -I + \gamma g_{\mathrm{h}} W_r$.
Continuous-time stability requires 
$\Re\!\left(\lambda_{\max}(J)\right) < 0$, and for Jacobians of this form
the associated instability boundary is well-approximated by the contour 
$r_e \approx 1$, 
where the reentry gain balances the unit leak.
The blue region therefore marks asymptotically stable trajectories,
while the yellow boundary identifies the onset of oscillatory or weakly
divergent flow in which reentry outweighs homeostatic regulation.

Together, these plots reveal that the network stability does not depend solely on the
magnitude of recurrence, but on the coordinated balance between feedback gain,
homeostatic normalization, and the spectral structure of $W_r$. This trade-off
defines a controllable stability manifold that enables reentrant computation
without runaway amplification.

\begin{figure}[ht!]
\includegraphics[scale=0.78, trim= 0.7cm 19cm 0cm 0cm]{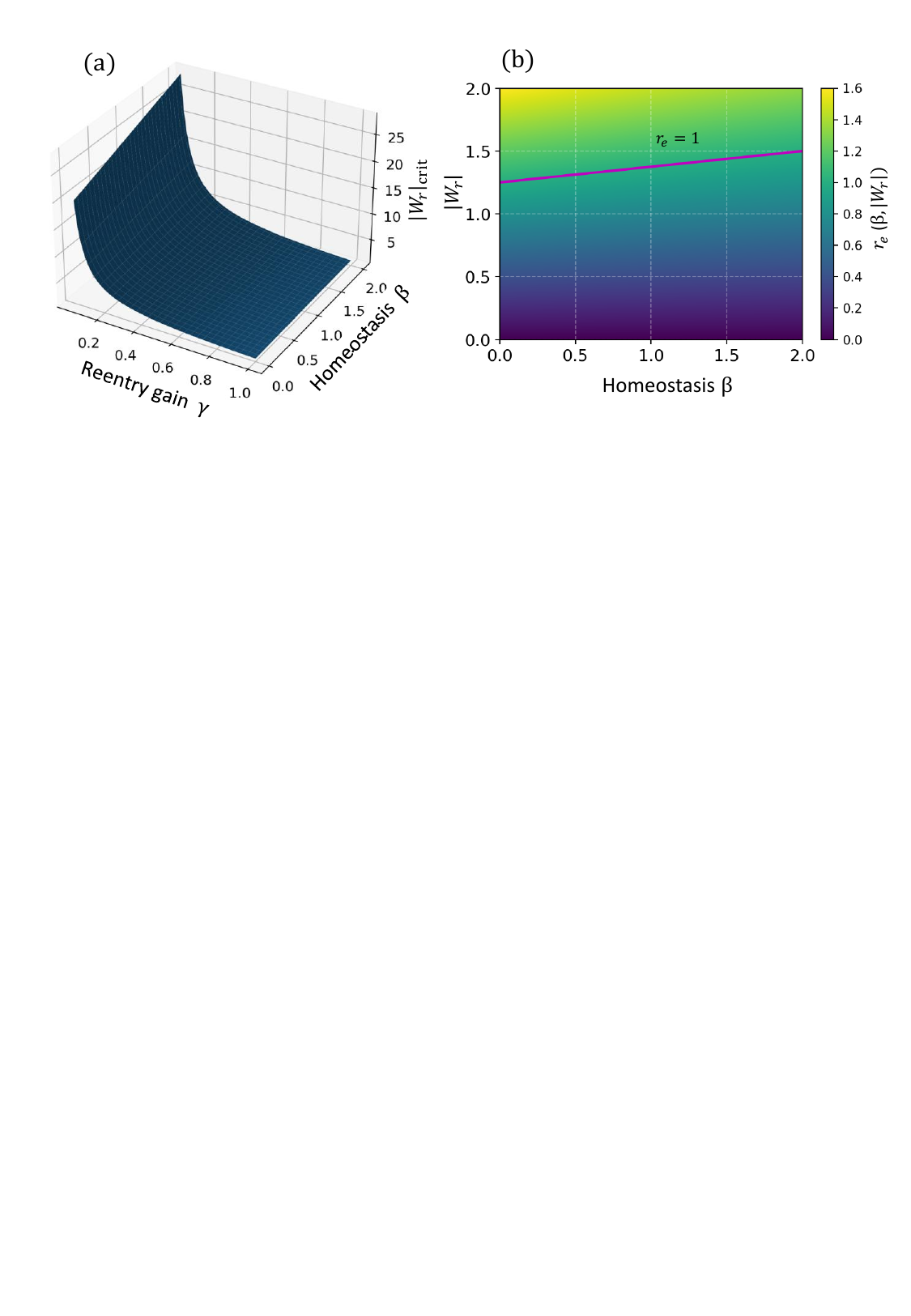}
\caption{
Stability structure of FHRN dynamics across the parameter space. 
(a) Critical stability surface showing the boundary 
$\|W_r\|_{\mathrm{crit}} = 1 / (\gamma\, g(\|y\|))$ across the 
$(\gamma,\beta)$ parameter space. Regions below the surface correspond 
to stable reentrant dynamics, while values above the boundary predict 
instability.  
(b) Two–dimensional stability map for fixed radius 
$\|y\|=1.1$, illustrating how stability varies with homeostasis $\beta$ 
and reentry magnitude $\|W_r\|$. Blue regions indicate contractive 
dynamics, while the transition curve ($r_e=1$) marks the 
onset of oscillatory or weakly divergent behavior.
}
\end{figure}

\section{Discussions and conclusion}

The continuous-time formulation developed here shows that the  \emph{Fast-Weights Homeostatic Reentry Network} occupies a computational regime not present in prior neural ODE frameworks.
Whereas conventional continuous recurrent systems propagate state through a fixed operator \citep{7,8}, 
and liquid neural networks stabilize dynamics via neuron-local adaptive time constants \citep{9}, 
FHRN introduces a state-dependent reentrant operator whose effective gain is regulated at the population level. 
This creates a system in which internal representations are not merely propagated but continuously revisited and reinterpreted, 
transforming recurrence from a temporal buffer into a mechanism for self-referential computation.

The resulting flow exhibits a bounded attractor geometry where rotational reentry, fast associative coupling, 
and radial homeostasis jointly determine the trajectory of population activity. 
This structure provides a mathematically grounded analogue of biological gain control, 
in which firing-rate homeostasis prevents runaway amplification while allowing oscillatory activity to persist \citep{22,23}. 
Unlike classical fast-weight models or attention-based retention mechanisms—which store associations 
but do not constrain how recursive inference unfolds—the reentry network tightly couples memory, feedback, 
and stability into a single, regulated dynamical system.

A notable consequence of this formulation is that recursive inference is not externally imposed but arises from the system’s own dynamics. 
As the reentry operator repeatedly encounters its evolving latent state, the network performs iterative refinement rather than producing a single deterministic mapping. 
This introduces the possibility of emergent reflective computation, 
wherein multi-step internal processing yields qualitatively distinct behavior as model scale, memory capacity, or reentry gain increase—analogous to emergent reasoning phenomena observed 
in large attention-based models \citep{24,25}.

Viewed through this lens, FHRN is not merely a parametric modification of fast-weight networks or a continuous analog of transformer retention layers. 
Instead, it represents a distinct computational principle: a dynamical system capable of maintaining and re-evaluating internal hypotheses through stable, self-referential flow. 
This framework suggests a path toward architectures that can sustain working memory, perform iterative internal inference, 
and potentially exhibit reflection-like dynamics—bridging concepts from biological recurrent computation and emergent behavior in large-scale artificial networks.

\emph{In conclusion},
the \emph{Fast-Homeostatic Reentrant Network} provides a unified dynamical formulation that integrates rapid feedback, short-term memory, and stability within a single continuous model. 
By introducing a population-level reentry operator and homeostatic gain, the framework transforms ordinary recurrent flows into \emph{self-referential yet bounded attractor dynamics}. 
This regime generalizes both fast-weight associative systems and liquid neural networks, while remaining analytically accessible through Jacobian and Lyapunov analysis. 
The resulting dynamics offer a compact physical model for cortical-like reentry: a mechanism that allows neural populations to maintain reflective internal states without instability. 
Beyond its theoretical implications, FHRN establishes a foundation for \emph{reflective computation}—machines that not only propagate information 
but dynamically reinterpret their own states in real time.

\section*{Acknowledgements}
This work was partially supported by the Institute of Information \& Communications Technology Planning \& Evaluation (IITP) grant funded by the Korea government (MSIT) (IITP-RS-2025-02214780).

The author acknowledges the support of ChatGPT (GPT-5, OpenAI) for assistance in literature review and conceptual structuring during early development.

\clearpage
\bibliographystyle{iclr2025_conference}

\begin{thebibliography}{99}

\bibitem[Abbott \& Regehr(2004)]{1}
Abbott, L. F. \& Regehr, W. G.
Synaptic computation.
\emph{Nature}, 431:796--803, 2004.

\bibitem[Gerstner et~al.(2014)]{2}
Gerstner, W., Kistler, W. M., Naud, R., \& Paninski, L.
\emph{Neuronal Dynamics}.
Cambridge University Press, 2014.

\bibitem[Turrigiano \& Nelson(2000)]{3}
Turrigiano, G. G. \& Nelson, S. B.
Hebb and homeostasis in neuronal plasticity.
\emph{Current Opinion in Neurobiology}, 10:358--364, 2000.

\bibitem[Mongillo et~al.(2008)]{4}
Mongillo, G., Barak, O., \& Tsodyks, M.
Synaptic theory of working memory.
\emph{Science}, 319:1543--1546, 2008.

\bibitem[Hopfield(1982)]{5}
Hopfield, J. J.
Neural networks and physical systems with emergent collective computational abilities.
\emph{Proceedings of the National Academy of Sciences}, 79:2554--2558, 1982.

\bibitem[Amit \& Gutfreund(1985)]{6}
Amit, D. J. \& Gutfreund, H.
Storing infinite numbers of patterns in a spin-glass model of neural networks.
\emph{Physical Review Letters}, 55:1530--1533, 1985.

\bibitem[Funahashi \& Nakamura(1993)]{7}
Funahashi, K. \& Nakamura, Y.
Approximation of dynamical systems by continuous time recurrent neural networks.
\emph{Neural Networks}, 6:801--806, 1993.

\bibitem[Beer(1995)]{8}
Beer, R. D.
On the dynamics of small continuous-time recurrent neural networks.
\emph{Adaptive Behavior}, 3:469--509, 1995.

\bibitem[Hasani et~al.(2021)]{9}
Hasani, R., Lechner, M., Amini, A., Rus, D., \& Grosu, R.
Liquid time-constant networks.
In \emph{Proceedings of AAAI Conference on Artificial Intelligence}, pp. 7657--7666, 2021.

\bibitem[Barak \& Tsodyks(2014)]{10}
Barak, O. \& Tsodyks, M.
Working models of working memory.
\emph{Current Opinion in Neurobiology}, 25:20--24, 2014.

\bibitem[Seeholzer et~al.(2019)]{11}
Seeholzer, A., Deger, M., \& Gerstner, W.
Stability of working memory in continuous attractor networks under the control of short-term plasticity.
\emph{PLoS Computational Biology}, 15, 2019.

\bibitem[Edelman(1989)]{12}
Edelman, G. M.
\emph{Neural Darwinism: The Theory of Neuronal Group Selection}.
Basic Books, 1989.

\bibitem[Tononi et~al.(1994)]{13}
Tononi, G., Sporns, O., \& Edelman, G. M.
A measure for brain complexity: Relating functional segregation and integration in the nervous system.
\emph{Proceedings of the National Academy of Sciences}, 91:5033--5037, 1994.

\bibitem[Hinton \& Plaut(1987)]{14}
Hinton, G. E. \& Plaut, D. C.
Using fast weights to deblur old memories.
In \emph{Proceedings of the 9th Annual Conference of the Cognitive Science Society}, pp. 177--186, 1987.

\bibitem[Schmidhuber(1992)]{15}
Schmidhuber, J.
Learning to control fast-weight memories: An alternative to dynamic recurrent networks.
\emph{Neural Computation}, 4:131--139, 1992.

\bibitem[Ba et~al.(2016)]{16}
Ba, J., Hinton, G. E., Mnih, V., Leibo, J., \& Ionescu, C.
Using fast weights to attend to the recent past.
arXiv:1610.06258, 2016.

\bibitem[Schlag et~al.(2021)]{17}
Schlag, I., Irie, K., \& Schmidhuber, J.
Linear Transformers are secretly fast weight programmers.
arXiv:2102.11174, 2021.

\bibitem[Katharopoulos et~al.(2020)]{18}
Katharopoulos, A., Vyas, A., \& Fleuret, F.
Transformers are RNNs: Fast autoregressive transformers with linear attention.
In \emph{Proceedings of NeurIPS}, 2020.

\bibitem[Sun et~al.(2023)]{19}
Sun, Y., Dong, L., Huang, S., Ma, S., Xia, Y., Xue, J., Wang, J., \& Wei, F.
Retentive Network: A successor to Transformer for large language models.
In \emph{Proceedings of ICLR}, 2023.

\bibitem[Irie et~al.(2022)]{20}
Irie, K., Schlag, I., Csordas, R., \& Schmidhuber, J.
A modern self-referential weight matrix that learns to modify itself.
arXiv:2202.05780, 2022.

\bibitem[Chae(2025)]{21}
Chae, B. G.
Recursive dynamics in fast-weights homeostatic reentry networks: Toward reflective intelligence.
arXiv:2511.06798, 2025.

\bibitem[Cannon \& Miller(2016)]{22}
Cannon, P. \& Miller, J.
Synaptic and intrinsic homeostasis cooperate to optimize single neuron response properties and tune integrator circuits.
\emph{Journal of Neurophysiology}, 116:2004--2022, 2016.

\bibitem[Niemeyer et~al.(2021)]{23}
Niemeyer, N., Schleimer, J. H., \& Schreiber, S.
Biophysical models of intrinsic homeostasis: Firing rates and beyond.
\emph{Current Opinion in Neurobiology}, 70:81--88, 2021.

\bibitem[Wei et~al.(2022)]{24}
Wei, J., Tay, Y., Bommasani, R., Raffel, C., Zoph, B., Borgeaud, S., Yogatama, D., Bosma, M., Zhou, D., Metzler, D., Chi, E. H., Hashimoto, T., Vinyals, O., Liang, P., Dean, J., \& Fedus, W.
Emergent abilities of large language models.
\emph{Transactions on Machine Learning Research}, 2022.

\bibitem[Webb et~al.(2023)]{25}
Webb, T., Holyoak, K. J., \& Lu, H.
Emergent analogical reasoning in large language models.
\emph{Nature Human Behaviour}, 7:1526--1541, 2023.

\end{thebibliography}

\clearpage
\appendix

\renewcommand{\theequation}{S\arabic{equation}}
\setcounter{equation}{0}

\section*{Appendix A: Lyapunov Stability Proof (Expanded)}

We analyze the stability of the continuous-time FHRN dynamics:
\begin{equation}
\dot{y} = -y + f(Wy, x, A) + g_{\mathrm{h}}(y),
\qquad 
\dot{A} = -\lambda A + \Phi(y,x),
\quad \lambda > 0,
\end{equation}
where $A$ denotes the fast-weight trace and $g_{\mathrm{h}}$ implements
population-level homeostatic regulation.

\medskip
We adopt the Lyapunov candidate:
\begin{equation}
V(y)=\frac{1}{4}(\|y\|^2-1)^2,
\end{equation}
which penalizes deviation of the activity norm from the equilibrium radius $\|y\|=1$.

\paragraph{(i) Time derivative.}
Using $\frac{d}{dt}\|y\|^2 = 2y^\top \dot{y}$, we obtain:
\begin{equation}
\dot{V} = (\|y\|^2 - 1)\, y^\top \dot{y}.
\end{equation}
Substituting the system dynamics yields:
\begin{equation}
\dot{V} 
= (\|y\|^2 - 1)\,y^\top[-y + f(Wy,x,A) + g_{\mathrm{h}}(y)].
\end{equation}

\paragraph{(ii) Leak contribution.}

The first component of the vector field arises from the intrinsic leak
term $-y$, which acts as a baseline stabilizing mechanism commonly found
in energy-based and continuous recurrent neural systems.  This term
continuously pulls the state toward the origin, counteracting unchecked
growth of the activity vector.

Its contribution to the Lyapunov derivative is straightforward to compute:
\[
y^\top(-y) = -\|y\|^2,
\]
and therefore
\begin{equation}
(\|y\|^2 - 1)\,y^\top(-y)
= -(\|y\|^2 - 1)\,\|y\|^2.
\end{equation}

Since $\|y\|^2 \ge 0$, we can bound this term more simply by observing that
$\|y\|^2 \ge 1$ whenever the radial deviation factor $(\|y\|^2-1)$ is positive.
Thus,
\begin{equation}
-(\|y\|^2 - 1)\|y\|^2
\;\le\;
-(\|y\|^2 - 1)^2,
\end{equation}
showing that the leak component contributes a strictly nonpositive term
to $\dot{V}$.  In other words, the leak introduces intrinsic dissipation:
whenever $\|y\|\neq 1$, it reduces the Lyapunov energy and therefore
promotes convergence toward the unit-norm manifold even in the absence
of homeostatic gain modulation or external drive.

\paragraph{(iii) Boundedness and regularity of the fast-weight drive.}

To ensure well-posedness of the continuous-time FHRN dynamics and
to enable Lyapunov-based stability analysis, we make a regularity
assumption on the nonlinear forcing term $f(Wy, x, A)$.
Specifically, we assume $f$ is \emph{Lipschitz sector-bounded} in the
state $y$, fast-weight trace $A$, and external drive $x$.
This corresponds to a standard stability constraint in nonlinear
recurrent systems and can be enforced in practice via normalization
mechanisms such as layer normalization, spectral norm control, or
activation clipping.

Formally, there exist nonnegative constants $L_y, L_A, L_x$ such that
\begin{equation}
\label{eq:lipschitz_f}
\|f(Wy,x,A) - f(0,0,0)\|
\;\le\;
L_y \|y\| + L_A \|A\| + L_x \|x\|.
\end{equation}

Without loss of generality we take $f(0,0,0)=0$ by centering the
representation space, yielding the simplified bound
\begin{equation}
\|f(Wy,x,A)\| \;\le\; L_y \|y\| + L_A \|A\| + L_x \|x\|.
\tag{\ref{eq:lipschitz_f}$'$}
\end{equation}

Applying the Cauchy--Schwarz inequality to the mixed term appearing
in the Lyapunov derivative,
\[
|y^\top f(Wy,x,A)| \le \|y\|\,\|f(Wy,x,A)\|,
\]
and substituting (\ref{eq:lipschitz_f}$'$) yields
\begin{equation}
\label{eq:sector_linear_form}
|y^\top f(Wy,x,A)|
\;\le\;
L_y \|y\|^2 + L_A \|y\|\,\|A\| + L_x \|y\|\,\|x\|.
\end{equation}

Since the FH-RL dynamics include a radial homeostatic feedback term
that regulates the magnitude of $y$, we additionally assume
$\|y\|\le R$ for some finite $R>0$ along system trajectories.
Under this boundedness condition, (\ref{eq:sector_linear_form})
admits the simplified linear sector bound
\begin{equation}
\label{eq:simplified_sector}
|y^\top f(Wy,x,A)|
\;\le\;
\alpha\|y\| + \beta_0\|A\| + \zeta\|x\|,
\end{equation}
where the constants are defined as aggregated Lipschitz parameters:
$\alpha := L_yR$, $\beta_0 := L_A R$, and $\zeta := L_x R$.

Finally, multiplying both sides of (\ref{eq:simplified_sector}) by
$|\|y\|^2-1|$ and applying Young's inequality produces bounded
quadratic terms compatible with the Lyapunov candidate
$V(y)=\tfrac14(\|y\|^2-1)^2$:
\begin{equation}
\label{eq:young_final_bound}
|(\|y\|^2 - 1)\,y^\top f(Wy,x,A)|
\;\le\;
c_f(\|y\|^2 - 1)^2 + c_A\|A\|^2 + c_x\|x\|^2,
\end{equation}
for suitably chosen nonnegative constants $c_f,c_A,c_x\ge 0$.

This bound ensures that the nonlinear forcing term can be treated as a
dissipative disturbance in the Lyapunov framework, completing the
conditions required for global asymptotic radial stability of the
FHRN dynamics.

\paragraph{(iv) Homeostatic damping.}

We now analyze the contribution of the homeostatic feedback term
\[
g_{\mathrm{h}}(y)= -\kappa(1-\|y\|^2)\,y,\qquad \kappa>0,
\]
which implements a radial gain control: when the population activity
is too small ($\|y\|<1$) the factor $(1-\|y\|^2)>0$ amplifies $y$,
whereas for overly large activity ($\|y\|>1$) the factor becomes
negative and suppresses $y$.  Thus $g_{\mathrm{h}}$ always acts to
restore the radius $\|y\|$ toward the homeostatic setpoint $\|y\|\approx 1$.

The inner product between the state and the homeostatic field is
\begin{equation}
y^\top g_{\mathrm{h}}(y)
= y^\top\big[-\kappa(1-\|y\|^2)y\big]
= -\kappa(1-\|y\|^2)\|y\|^2.
\end{equation}
This term is nonpositive for all $y$ and vanishes only on the
homeostatic manifold $\|y\|=1$.

In the Lyapunov derivative, $y^\top g_{\mathrm{h}}(y)$ appears
multiplied by the radial deviation factor $(\|y\|^2-1)$:
\begin{equation}
(\|y\|^2 - 1)\,y^\top g_{\mathrm{h}}(y)
= (\|y\|^2 - 1)\big[-\kappa(1-\|y\|^2)\|y\|^2\big]
= -\kappa(\|y\|^2 - 1)^2 \|y\|^2.
\end{equation}
Since $(\|y\|^2-1)^2\ge 0$ and $\|y\|^2\ge 0$, this contribution is
always nonpositive:
\begin{equation}
(\|y\|^2 - 1)\,y^\top g_{\mathrm{h}}(y)
\;\le\; 0,
\end{equation}
and, in particular, for $\|y\|\ge 1$ we obtain the simpler bound
\begin{equation}
(\|y\|^2 - 1)\,y^\top g_{\mathrm{h}}(y)
= -\kappa(\|y\|^2 - 1)^2 \|y\|^2
\;\le\; -\kappa(\|y\|^2 - 1)^2,
\end{equation}
because $\|y\|^2\ge 1$ in this regime.
Thus the homeostatic term provides a strictly dissipative
(radially contracting) contribution away from the unit-norm
manifold, directly counteracting any destabilizing fast-weight
drive in the Lyapunov analysis.

\paragraph{(v) Combined inequality.}

Substituting the bounds derived in (iii) and the dissipative
contribution of the homeostatic field from (iv) into the Lyapunov
derivative yields a competition between stabilizing and destabilizing
terms.  The fast-weight and external-drive terms introduce bounded,
potentially destabilizing components, while the homeostatic damping
contributes a negative quadratic term that contracts the dynamics
whenever the population activity deviates from its target radius.

Collecting all contributions, we obtain the aggregated inequality
\begin{equation}
\dot{V}(y)
\;\le\;
-(1+\kappa-c_f)\,(\|y\|^2 - 1)^2
\;+\;
c_A\|A\|^2
\;+\;
c_x\|x\|^2,
\label{eq:combined_lyapunov}
\end{equation}
where the coefficient $(1+\kappa-c_f)$ represents the net strength of
the intrinsic radial contraction relative to the perturbation induced
by recurrent and fast-weight feedback.

A sufficient condition for strict Lyapunov contraction is therefore
\begin{equation}
1+\kappa-c_f > 0
\quad\Longleftrightarrow\quad
\kappa > c_f - 1,
\end{equation}
ensuring that the dissipative homeostatic term dominates the Lipschitz-bounded
nonlinear contributions. Under this parameter regime, $\dot{V}<0$ whenever
$\|y\|\ne 1$, establishing global radial asymptotic stability of the
FHRN dynamics in the absence of external drive.

\paragraph{(vi) Stability of $A$ and ISS.}
Since
\begin{equation}
\dot{A} = -\lambda A + \Phi(y,x),\qquad \lambda>0,
\end{equation}
standard linear system theory implies $A$ remains bounded given bounded $y$ and $x$.
Thus $V(y)$ is bounded and trajectories converge to an invariant set satisfying $\|y\|\approx1$.
The system is \textit{input-to-state stable (ISS)} with respect to $(x,A)$.
Leak and homeostasis ensure radial damping, while reentry via $W_r$ contributes only rotational structure,
yielding a globally bounded reflective flow on the activation manifold.

\section*{Appendix B: From Oja’s Rule to Radial Homeostasis in FHRN}

To clarify the mathematical origin of the homeostatic term used in FHRN, 
we show that it can be derived directly from Oja’s normalized Hebbian learning rule via a change of variables from weight space to activity space. 
This establishes a formal correspondence between synaptic normalization and radial stabilization of neural population dynamics.

\paragraph{(i) Scalar derivation.}
Consider a single neuron with constant input $x$, weight $w$, and activation
\begin{equation}
y = wx.
\end{equation}
Oja’s rule for normalized Hebbian learning is:
\begin{equation}
\dot{w} = xy - y^2 w.
\label{eq:oja}
\end{equation}
Substituting $y = wx$ into Eq.~\eqref{eq:oja} yields:
\begin{equation}
\dot{w} = x(wx) - (wx)^2 w = x^2 w(1 - w^2).
\end{equation}
Since $y = wx$, we have $\dot{y} = x \dot{w}$, giving:
\begin{equation}
\dot{y} = x \cdot x^2 w(1 - w^2) = x^2 y - y^3 = y(x^2 - y^2).
\label{eq:y_scalar}
\end{equation}
Interpreting $x^2$ as a target activation radius $\theta^2$, Eq.~\eqref{eq:y_scalar} becomes:
\begin{equation}
\dot{y} = y(\theta^2 - y^2) = -y(y^2 - \theta^2).
\label{eq:scalar_homeostasis}
\end{equation}
Equation~\eqref{eq:scalar_homeostasis} represents a radial restoring force that stabilizes the activity magnitude at $\theta$.

\paragraph{(ii) Vector generalization.}

For a neural population state $y \in \mathbb{R}^d$, the scalar expression generalizes naturally by replacing $y^2$ with $\|y\|^2$:
\begin{equation}
\dot{y} = -\kappa\left(\|y\|^2 - \theta^2\right)y,
\label{eq:fh-homeo}
\end{equation}
where $\kappa > 0$ controls the homeostatic rate. Eq.~\eqref{eq:fh-homeo} is precisely the polynomial homeostatic regulator employed in the neural network 
to maintain bounded population activity under recurrent reentry.

\paragraph{(iii) Interpretation.}

Oja’s term $-y^2w$ is traditionally interpreted as a weight normalization mechanism preventing runaway strengthening of synapses. 
Under the substitution $y = Wx$ and fixed input assumptions, the same term becomes a radial inhibitory factor $-(\|y\|^2 - \theta^2)y$ governing the stability of neural activity, not the weights. 
Thus:
\begin{equation}
-y^2w \quad \longleftrightarrow \quad -(\|y\|^2 - \theta^2)y,
\end{equation}
revealing that weight normalization and activity homeostasis are mathematically equivalent processes expressed in different coordinate spaces.
This derivation demonstrates that FHRN’s homeostatic term arises naturally as the activity-space analogue of Oja’s normalized Hebbian learning rule. To our knowledge, this is the first explicit formulation showing that associative normalization in synaptic space implies radial stability in population-state dynamics.

\end{document}